\documentclass{article}
\title{$A$-Hypergeometric Distributions \\ and Newton Polytopes}
\author{Nobuki Takayama, Satoshi Kuriki, Akimichi Takemura}
\date{January 4, 2015, Revised November 10, 2015}
\usepackage{amssymb,amsmath,amsthm}
\usepackage[dvipdfmx]{hyperref}
\usepackage[dvipdfmx]{color}
\usepackage{ulem}  
\usepackage{array,arydshln}
\newcommand\Z{\mathbb{Z}}
\newcommand\R{\mathbb{R}} 
\newcommand\A{{\bar A}}
\newcommand\N{\mathbb{N}}
\renewcommand\a{{\bar a}}

\def\kappa{k}

\def\qqq#1{{\color{red} #1}}
\def\comment#1{ }
\def\pd#1{\partial_{#1}}
\def\qed{\hskip 15pt /\kern -2pt/ \bigbreak}
\newtheorem{example}{Example}
\newtheorem{proposition}{Proposition}
\newtheorem{theorem}{Theorem}

\newtheorem{lemma}{Lemma}
\newtheorem{corollary}{Corollary}

\begin{document}
\maketitle

Abstract:
We give a bijection between a quotient space of the parameters
and the space of moments for any $A$-hypergeometric distribution.
An algorithmic method to compute the inverse image of the map
is proposed utilizing the holonomic gradient method and
an asymptotic equivalence of the map and the iterative proportional scaling.
The algorithm gives a method for solving a conditional maximum likelihood estimation problem in statistics.
The interplay between the theory of hypergeometric functions and statistics allows us to 
give some new formulas for $A$-hypergeometric polynomials.

\section{Introduction}  \label{sec:intro}
We denote by $\N$ the set of the non-negative integers.
Let $A$ be a $d \times n$ 
configuration matrix with non-negative integer entries.
We assume that the rank of $A$ is $d$.
The  $A$-hypergeometric polynomial \cite{sst-c}
for $A$ and $\beta \in \N^d$
is defined by 
\begin{equation}
\label{eq:hypergeometric-polynomial}
 Z(\beta;p) = \sum_{ Au=\beta, u \in \N^n} \frac{p^u}{u!},
\end{equation}
where
$p^u = \prod_{i=1}^n p_i^{u_i}$ and
$u! = \prod_{i=1}^n u_i!$.
Set $p_i = \exp \xi_i$
and let $\exp \xi$ denote the vector $(\exp \xi_1, \ldots, \exp \xi_n)$.
We fix $\beta \not= 0$ such that $\beta \in \N A = \sum_{i=1}^n \N a_i$,
where $a_i$ denotes the $i$-th column vector of the matrix $A$.
Let $U \in \N^n$ be a random variable of the $(A,\beta)$
{\it hypergeometric distribution with the parameter} $p \in \R_{>0}^n$ (or $\xi\in \R^n$),
which is defined by 
\begin{equation}
\label{eq:a-distribution}
 P(U=u \,|\, Au=\beta ) = \frac{p(\xi)^u}{u! Z(\beta;p(\xi))}
 = \frac{\exp( u\cdot \xi)}{u! Z(\beta;p(\xi))},
\quad u\cdot \xi = \sum_{i=1}^n u_i \xi_i.
\end{equation}
If no confusion arises, we simply call this the $A$-hypergeometric distribution.
The $A$-hypergeometric distribution is in turn a generalization of
the generalized ($p_i$ may take any positive number) 
hypergeometric distribution on the contingency tables with fixed marginal
sums (see, e.g., \cite[Chapters 4, 6]{dojo}, \cite{ott}), 
and  is the conditional distribution of $u$ given by $\beta =Au$ under
the Poisson distribution
\begin{equation}
\label{eq:unconditional}
 P(U=u) = \frac{p^u}{u!}\exp(-{\bf 1} \cdot p), \ \ {\bf 1}=(1,\ldots,1).
\end{equation}
In the setting of testing statistical hypotheses, this corresponds to the alternative hypothesis against the null hypothesis of $p_i=s^{a_i}$, $s \in \R_{>0}^d$, $\forall i$. The polynomial $Z$ is the {\it normalizing constant} or the {\it partition function} of the $A$-hypergeometric distribution.

The expectation of the random variable $U_i$ under \eqref{eq:a-distribution}
is equal to
\begin{equation}  \label{eq:expectation1}
 \sum_{Au=\beta, u \in \N^n} u_i \frac{p(\xi)^u}{u! Z(\beta;p(\xi))}. 
\end{equation}
Setting 
$$\psi(\xi) = \log Z(\beta;p(\xi)), $$
the expectation of $U_i$ is written as
$$ E[U_i] = 
 \frac{p_i \pd{i} \bullet Z}{Z} |_{p=p(\xi)}
= \frac{\partial \log \psi(\xi)}{\partial \xi_i} ,
$$
where $p = (\exp \xi_1, \ldots, \exp \xi_n)$ and 
$\pd{i} = \frac{\partial}{\partial p_i}$.
If we set $\eta_i = E[U_i]$ and $\eta=(\eta_i)$,
which is a function of $\xi$, then the $\xi$-space and the $\eta$-space are dual
by the {\it moment map} $E[U]$
in the context of
the information geometry \cite{amari}.

We study here the map
between the $\xi$-space (the space of the parameters)
and the $\eta$-space (the space of moments).
This correspondence has been studied from several points of view in statistics and information geometry, e.g., 
\cite{amari}, \cite{barndorff-nielsen}, \cite{brown}, \cite{fienberg-rinaldo}, \cite{haberman}.
In section \ref{sec:expectation},
we determine the image of the $\xi$-space $\R^n$
by the moment map in the $\eta$-space,
which is described in terms of the Newton polytope of the polynomial $Z$.
We introduce a quotient space, which is called the space of
{\it the generalized odds ratios},  
of the $\xi$-space and construct an isomorphism
between the quotient space and the Newton polytope in the $\eta$-space
in Theorem \ref{th:one-to-one}.
This also yields, as a byproduct, a new theorem
on the image of the logarithmic derivatives of the Lauricella hypergeometric polynomial $F_D$ 
which is given in Theorem \ref{th:fdimage} in  section \ref{sec:classical}.
In section \ref{sec:mle}, we discuss an iteration method for computing
the inverse map of the moment map from the $\eta$-space to the $\xi$-space
with the help of the holonomic gradient method (HGM) \cite{N3OST2}, \cite{ot2},
which amounts to obtaining the conditional maximal likelihood estimate (MLE).
A subclass of this problem for $2\times m$ contingency tables is discussed in \cite{ott}.
We also note that several important quantities
in information geometry
such as Fisher matrices
can be numerically evaluated by the HGM in the case of
$A$-hypergeometric distributions.
In section \ref{sec:IPS}, we discuss the generalized odds ratio and a log-affine model to give the iterative proportional scaling (IPS) method
in Theorem \ref{th:ips}.
Finally, in section \ref{sec:equivalence}, we present a theorem for approximating
the normalizing constant and the inverse map by the IPS method, which is used as an initial value of the iteration
proposed in section \ref{sec:mle}.
We also note that this approximation theorem (Theorem \ref{th:EisM}) gives a new approximation
formula (Theorem \ref{th:approxZ}) for $A$-hypergeometric polynomials and consequently
for hypergeometric polynomials in several variables.

\section{Expectation and Newton Polytope}  \label{sec:expectation}

We are interested in the image of the map $E[U]=(E[U_1], \ldots, E[U_n])$
as a function of $\xi \in {\R}^n$
$$ E[U] \ :\ {\R}^n \ni \xi \mapsto \eta=E[U](\xi) \in {\R}^n $$
or a function of $p \in {\R}_{>0}^n$
$$ E[U] \ :\ {\R}_{>0}^n \ni p \mapsto \eta=E[U](p) \in {\R}^n. $$
We denote both functions by the same symbol $E[U]$ 
as long as no confusion arises.

\begin{proposition}  \label{prop:image1}
When $\eta$ is an image of $\xi \in {\R}^n$ 
by the moment map $E[U]$, 
we have $A \eta = \beta$ and $\eta \in {\R}_{\ge 0}$.
\end{proposition}

{\it Proof}\/.
This is an immediate consequence of (\ref{eq:expectation1}), because each $u$ in
the summand of (\ref{eq:expectation1}) satisfies $Au=\beta$.
\qed

We will call the polytope defined by $A \eta=\beta$, $\eta \in {\R}_{\geq0}$
a {\it generalized transportation polytope}\/.
We call its relative interior an {\it open generalized transportation polytope}\/.
As long as no confusion arises, we simply call the generalized transportation polytope
the {\it transportation polytope}\/.

We denote by $a_i \in \N^d$ the $i$-th column vector of $A$.
\begin{proposition} \label{prop:invariance}
\begin{enumerate}
\item 
The $A$-hypergeometric distribution (\ref{eq:a-distribution}) as a function of $p$ is invariant under the torus action  of $A$. In particular, 
the moment map $E[U]$ as a function of $p$ is invariant under the torus action of $A$.
In other words, fixing a vector  $y>0$, 
$E[U](y_1 s^{a_1}, \ldots, y_n s^{a_n})$ is a constant for all $s \in {\R}_{>0}^d$.
\item
The $A$-hypergeometric distribution (\ref{eq:a-distribution}) and the
 moment map $E[U]$ as a function of $\xi$
are constant on the image of $A^T\,:\, {\R}^d \rightarrow {\R}^n$.
In particular,
when $\xi - \xi' \in {\rm Im}\, A^T$,
$E[U](\xi)=E[U](\xi')$.
\end{enumerate}
\end{proposition}

{\it Proof}\/.
Substitute $p$ by $(y_1 s^{a_1}, \ldots, y_n s^{a_n})$ in 
(\ref{eq:a-distribution}).
Then, we obtain item 1  of the proposition.
Setting $y_i = \exp(\xi_i)$ and $s_i = \exp(\sigma_i)$ in (1),
we have $y_i s^{a_i} = \exp(\xi_i + a_i \cdot \sigma)$,
where $\sigma = (\sigma_1, \ldots, \sigma_d)$.
Since
$(\sigma \cdot a_1, \ldots, \sigma \cdot a_n) \in {\rm Im}\, A^T$,
we have item 2.
\qed

Proposition \ref{prop:invariance} implies that the function $E[U](\xi)$
can be defined on $\R^n / {\rm Im}\, A^T$.
Let us give a description of this quotient space.
We regard the matrix $A$ as a map from
$\Z^n$ to $\Z^d$ of maximal rank.
By suitable unimodular matrices $R$ on $\Z^n$ and $S$ on $\Z^d$,
we have the Smith normal form
\def\hsymb#1{\mbox{\strut\rlap{\smash{\Huge$#1$}}\quad}}
$$
S A R =
\left(
\begin{array}{ccc:c}
 \alpha_1 &         & O        &  \\
          & \ddots  &          & \ O \\
 O        &         & \alpha_d &  \\  
\end{array}
\right),\ \ %
 \alpha_i \not= 0,\ \alpha_i | \alpha_{i+1}.
$$
A $\Z$-module basis of ${\rm Ker}(A : \Z^n \rightarrow \Z^d)$ is
$\{ R e_{d+1}, \ldots, R e_n\}$,
where $\{e_i\}$ is the standard basis of $\Z^n$
expressed as column vectors.
We denote the vector $(R e_{d+i})^T$ by 
${\bar a}_i$.
Define a matrix $\A$ as
\[
 \A = \begin{pmatrix} \a_1 \\ \vdots \\ \a_{n-d} \end{pmatrix}_{(n-d)\times n},
\]
which is called the Gale transform of $A$.
Setting $\lambda = \A \xi$,
since 
$\{ {\bar a}_i^T \,|\, i=1, \ldots, n-d \}$
is also a basis of ${\rm Ker}(A : \R^n \rightarrow \R^d)$,
the map 
$$ \R^n / {\rm Im}\, A^T \ni \xi \mapsto \lambda \in \R^{n-d} $$
is an isomorphism.
We call the ratio of $p_j$'s
$\exp(\lambda_i) = p^{{\bar a}_i} = \prod_{j=1}^n p_j^{{\bar a}_{ij}}$
{\it the generalized odds ratio}.
We will discuss this ratio 
in section \ref{sec:IPS}.

When $E[U]$ is a function in one variable $t$ modulo ${\rm Im}\, A^T$,
the image can be determined in an elementary way.
In fact,
if we set $F(t) = \sum_{k=0}^N c_k \exp(kt)$, where $c_k \geq 0$,
then we have the following lemma.
\begin{lemma}  \label{lm:exppoly}
The function $\frac{d}{dt} \log e^{Mt} F(t)$ is an increasing function on ${\R}$ for any real number $M$.
\end{lemma}

{\it Proof}\/.
$$\frac{d^2}{dt^2} \log e^{Mt} F(t)
= \frac{F'' F - (F')^2}{F^2}. $$
The numerator is
$$ \left( \sum_{i=0}^N c_i i^2 \exp(it) \right)
   \left( \sum_{j=0}^N c_j \exp(jt) \right)
 - 
   \left( \sum_{i=0}^N c_i i \exp(it) \right)
   \left( \sum_{j=0}^N c_j j \exp(jt) \right).
$$
Expanding the products, the coefficient of $\exp(lt)$, $0\le l\le 2N$, is equal to
$$
   \sum_{i+j = l} (i^2 c_i c_j - i j c_i c_j)
 = \sum_{i+j=l, i>j} (i^2 - i j + j^2 - ij) c_i c_j 
 = \sum_{i+j=l, i>j} (i-j)^2 c_i c_j \geq 0
 $$
Therefore, $\frac{d}{dt} \log e^{Mt} F(t)$ is an increasing function.
\qed

\begin{example}\rm  \label{ex:2F1a}
$$
A = \left(
\begin{array}{cccc}
1 & 1 & 0 & 0 \\
0 & 0 & 1 & 1 \\
1 & 0 & 1 & 0 \\
\end{array}
\right), \quad
\beta=(b_1,b_2,c_1).
$$
Here, $u$ satisfying $Au=\beta$ can be regarded as a $2 \times 2$ contingency table
$$ u =
\left(
\begin{array}{cc}
 u_1 & u_2 \\
 u_3 & u_4 \\
\end{array}
\right)
=
\left(
\begin{array}{cc}
 u_{11} & u_{12} \\
 u_{21} & u_{22} \\
\end{array}
\right)
$$
with row sums $b_1$ and $b_2$
and column sums $c_1$ and $c_2 = b_1 + b_2 -c_1$.
The invariance under the torus action and $\A = (-1,1,1,-1)$ imply that 
$E[U](p)$ depends only on
the odds ratio $z=\frac{p_{12} p_{21}}{p_{11} p_{22}}$, where $p_1 = p_{11}, p_2=p_{12}, p_3=p_{21}$, and $p_4=p_{22}$.

Let us illustrate how the expectation is expressed in terms of a hypergeometric series
by an example.
When $(b_1,b_2,c_1,c_2) = (36,12,37,11)$,
the expectation $E[U_{11}](p)$ is equal to
$$ 36 - \frac{z F'(z)}{F(z)}, \quad
 F(z) = F(-36,-11,2; z).
$$
Here, $F(z)$ is the Gauss hypergeometric polynomial or can be regarded as a Jacobi polynomial.
In fact,
we have
$$ Z(p) = \frac{p^\mu}{\mu!} 
 F\left(-b_1,-c_2,c_1-b_1+1; \frac{p_{12}p_{21}}{p_{11} p_{22}}\right), \quad
\mu = \left(
\begin{array}{cc}
 -b_1    & 0 \\
 c_1-b_1 & -c_2 \\
\end{array}
\right).
$$
The expectation 
$E[U_{11}] = p_{11} \frac{\partial Z}{\partial p_{11}} / Z$
is equal to
$$
 \mu_{11} 
+ \frac{p^{\mu}}{\mu!} \left(-\frac{p_{12} p_{21}}{p_{11} p_{22}}\right) F'(z)/Z 
 = \mu_{11}-z \frac{d F/dz}{F}.
$$

Let us determine the image of the moment map $E[U]$.
It follows from Proposition \ref{prop:image1} that the image lies on the domain
$ 25 < E[U_{11}] < 36$.
Since $\log \psi(\xi)$ is a lower convex function 
by Lemma \ref{lm:exppoly} or
by a general theorem for the exponential family
and
$ \frac{\partial E[U_{11}](z(\xi))}{\partial \xi_{11}} =
 - \exp(\xi_{12}+\xi_{21}-\xi_{11}-\xi_{22}) \frac{d E}{dz}(z(\xi))$,
we have $\frac{d E}{dz} \leq 0$. 
Therefore, $E[U_{11}](z)$ is a decreasing function in $z$.
We have $F(z) = 1 + \cdots + \frac{(-36)_{11} (-11)_{11}}{(2)_{11} 11!} z^{11}$.
Taking the limit $z \rightarrow 0$, we have
$\frac{z F'(z)}{F(z)} \rightarrow 0$
and taking the limit $z \rightarrow +\infty$, we have
$\frac{z F'(z)}{F(z)} \rightarrow 11$. Then, the expectation converges to $36 - 11=25$. 
Thus, the image agrees with the interval $(25,36)$ and the vertices $36$ and $25$ are attained by
$p = [[1,0],[1,1]]$ and $p = [[0,1],[1,1]]$ respectively.
In other words, the image for $p \in {\R}_{\geq 0}^4$ is $[25,36]$.

Analogously, we have $E[U_{12}]=z F'/F$ and the image is
$(0, 11)$. 
\end{example}

This example in the one-variable case can be generalized as follows.

\begin{theorem}  \label{th:one-to-one}
\begin{enumerate}
\item
The image of the moment map $E[U]$ agrees with the relative interior of 
the Newton polytope ${\rm New}(Z)$ of the normalizing constant $Z$ as a polynomial in $p$
when the dimension of the Newton polytope is $n-d$.
\item
The map 
$$E[U] \ : \ {\R}^n/{\rm Im}\,A^T \longrightarrow {\rm relint}({\rm New}(Z)) $$
is one-to-one
when the dimension of the Newton polytope is $n-d$,
where ``relint'' denotes the relative interior.
\end{enumerate}
\end{theorem}
Before proceeding to the proof, 
we note two sufficient conditions so that the dimension of the Newton polytope
is $n-d$.
Let $q_i$ be the non-negative integer 
${\rm max}_{j=1, \ldots, n-d} |({\bar a}_j)_i|$,
where $({\bar a}_j)_i$ denotes the $i$-th entry of the vector ${\bar a}_j$.
If $\beta \in \N A = \sum_{i=1}^n \N a_i$ lies in
\begin{equation} \label{eq:dimcond}
 \R_{\geq q} A = \sum_{i=1}^n \R_{\geq q_i} a_i, \quad
\R_{\geq q_i} = \{ c \in \R \,|\, c \geq q_i \},
\end{equation} 
then the dimension of the Newton polytope is $n-d$.

The second sufficient condition is $\beta = k \beta'$
for $k$ being a sufficiently large natural number 
and $\beta' \in \N A$ satisfying $\beta' \in {\rm int}(\R_{\geq 0} A)$.
This condition follows from condition 1 and $\beta'_i \not = 0$
for all $i$.

{\it Proof of Theorem \ref{th:one-to-one}}\/.
We regard $Z$ as a polynomial in $p$. Since the case of a monomial $Z$ is trivial,
we consider the case that $Z$ is not a monomial.
We denote by $S(Z)$ the support of $Z$. 
We will prove that when $\eta$ is in the relative interior of 
the Newton polytope ${\rm New}(Z)$,
which is the convex hull of $S(Z)$,
there exists an inverse image of $\eta$ by the moment map.
Let $m$ be a vertex of the Newton polytope ${\rm New}(Z)$.
We note that it is contained in the (closed) transportation polytope.
If we set $f = p^m/m!$ and $g = Z-f$,
then $\log Z = \log f (1 + \frac{g}{f})$.
Since the coefficients of the expansion of $Z$ are positive,
we have $\frac{g}{f} > 0$.
Therefore, we have
$$\log Z \geq \log f= m \cdot \log p -\log m!. $$
Then, $-\log Z$ is bounded as
$$ - \log Z(\beta;\xi) \leq - m \cdot \xi +\log m!. $$
Let $\eta$ be a point of the open Newton polytope of $Z$,
by which we mean the relative interior of the Newton polytope.
We consider the cost function 
\begin{equation} \label{eq:loglike}
 f(\xi)=\eta \cdot \xi - \log Z(\beta,p(\xi)). 
\end{equation}
The partial derivative $\frac{\partial f}{\partial \xi_i}$
is equal to $\eta_i - \frac{\partial Z}{\partial \xi_i}/Z$.
Then, $E[U_i]=\eta_i$ is equivalent to 
the partial derivative being equal to $0$.
Then, the existence of the maximum of $f(\xi)$ at the point $\xi=q$
implies that
${\rm grad}(f)=0$ holds at point $q$.
We will prove that the cost function has a maximum.
The cost function is bounded above by
$ \eta \cdot \xi - m \cdot \xi + \log m!$.
Let $C_m$ be the outer normal cone of ${\rm New}(Z)$.
In other words,
$C_m = \{ w \in {\R}^n \,|\, (y - m) \cdot w \leq 0 
 \mbox{ for all $y \in {\rm New}(Z)$} \}$.
Since the dimension of the Newton polytope is $n-d$,
the cone $C_m$ contains the linear space ${\rm Im}\, A^T$
which is orthogonal to the elements of the kernel of $A$
and is maximal.
Let $\theta$ be the angle between the two vectors $\xi$ and $\eta-m$.
When $\xi \in C_m$ and $\xi \not\in {\rm Im}\, A^T$, 
we have $\cos \theta <0$.
Since $\eta$ is a point in the interior of ${\rm conv}(S(Z))$,
there exists $\varepsilon_m$ such that
$ \cos \theta < \varepsilon_m < 0$ holds for any $\xi \in C_m$, $\xi \not\in {\rm Im}\, A^T$.
Let $M$ be a negative number which is smaller than
${\rm sup}\, (\eta \cdot \xi - \log Z)$.
Then, if $|\xi| > \frac{M}{\varepsilon_m | \eta-m |}$
and $\xi \not \in {\rm Im}\,A^T$,
we have 
$$M > |\xi| \varepsilon_m |\eta-m| > (\eta-m) \cdot \xi. $$
Let $\varepsilon$ be the maximum of $\varepsilon_m$ over all vertices $m$
of ${\rm conv}(S(Z))$.
Then, for any $\xi \in {\R}^n$, $\xi \not\in {\rm Im}\, A^T$,
the condition
$|\xi| > \frac{M}{\varepsilon} {\rm max} |\eta-m| $
implies $\eta \cdot \xi - \log Z < M$.
This means that when $E[U](\xi)$ is regarded as a function on ${\R}^n/{\rm Im}\, A^T$,
the function value is smaller than $M$ outside a compact domain.
Therefore, the function $\eta \cdot \xi - \log Z$ has a maximum.
Let 
$$ \xi(\eta) = {\rm maxarg}_\xi \, \left( \eta\cdot \xi - \log Z(\beta,p(\xi)) \right). $$ 
At the point $\xi(\eta)$, we have
$\eta = {\rm grad}(\log Z)$,
because the partial derivatives of $\log Z$ vanish at $\xi$ 
by the maxarg property.
We have proved that the inverse image of $\eta \in {\rm relint}({\rm New}(Z))$ 
exists.

Let us show that the image lies in the relative interior of ${\rm New}(Z)$.
Let $\eta$ be on the complement of ${\rm New}(Z)$.
There exists a facet hyperplane $L$ such that
$L(\eta) < 0$ and any point $u$ in ${\rm New}(Z)$ lies in the opposite
side of $\eta$, in other words $L(u) \geq 0$ holds.
Take a vector $m$ on $L=0$,
and let $C_m$ be the outer normal cone of $m$.
For any $\xi \in C_m$ and $u \in {\rm New}(Z)$, 
$ \xi \cdot (u-m) \leq 0$.
Therefore, ${\rm in}_\xi(Z)$ contains the term $p^m$.
Take $\xi^1 \in C_m$ such that $\xi^1 \cdot (\eta-m) > 0$,
and consider the cost function 
$f(\xi)=\xi \cdot \eta - \log Z = \eta \cdot \xi - m \cdot \xi - \log(Z/p^m)$.
Let $t$ be a scalar variable.
We restrict the cost function to the one-dimensional vector space
parameterized as $\xi=\xi^1 t$.
If no confusion arises, we denote by $f(t)$ the restricted cost function.
Then, we have 
$$f(t) = t\xi^1 \cdot (\eta-m)
 - \log \sum_u \exp(t \xi^1 \cdot (u-m))/u!.
$$
Since $\xi^1 \in C_m$, we have $\xi^1 \cdot (u-m) \leq 0$.
On the other hand, we have $\xi^1 \cdot (\eta-m) >0$.
Therefore, $f(t) \rightarrow +\infty$ when
$t \rightarrow +\infty$,
which means that cost function $f(\xi)$ does not have a maximum.
Since $f(\xi)$ is upper convex and smooth, this implies that
${\rm grad}(f)$ is not the zero vector at any point (see Lemma \ref{lm:zero} below).

Finally, we consider the case when $\eta$ is on the boundary of ${\rm New}(Z)$.
We suppose that $L(\eta)=0$ and that 
$L(u) \geq 0$ holds for any point $u$ of ${\rm New}(Z)$.
Let $m$ be a vertex of the Newton polytope on the hyperplane $L(e)=0$.
Let $\xi^1$ be a vector on the border of the outer normal cone $C_m$
such that it is orthogonal to the hyperplane $L(e)=0$.
Then, we have $\xi^1 \cdot (\eta-m)=0$.
We suppose that $f(\xi)$ has a maximum at $\xi=q$.
We restrict it to $\xi = q + t \xi^1$ and denote by $f(t)$
the restricted cost function.
We have
$f (t) = (q+t \xi^1) \cdot \eta - m \cdot(q+t\xi^1)-\log Z/p^m
 = q \cdot (\eta-m)-\log Z/p^m$.
Since $Z/p^m = \sum_u \exp(u\cdot q + t u \cdot \xi^1 - m \cdot q - t m \cdot \xi^1)/u!
= Z/p^m = \sum_u \exp((u-m)\cdot q) \exp(t (u -m) \cdot \xi^1)/u!
$,
$|f(t)|$ is bounded when $t \rightarrow +\infty$.
We note that $f(0)$ is the maximum from the assumption.
Then, we have $f(0) \geq f(t)$ for all $t$.
Since the terms in $Z/p^m$ are positive, $f(t)$ is not a constant function.
Since  $f(t)$ is holomorphic
and upper convex, there exists $t^0 >0$ such that $f(0)>f(t^0)$.
Therefore, there exists $t^1$ such that
$f'(t^1) = (f(t)-f(0))/t < 0$.
Since $f(t)$ is upper convex, $f''(t) \geq 0$ and consequently
$f'$ is not increasing function.
Therefore, for $t \geq t^1$, we have $f'(t) \leq f'(t^1) < 0$.
This implies that $f(t) \rightarrow -\infty$ when $t \rightarrow +\infty$.
This contradicts the assumption that $|f(t)|$ is bounded.

Let us show 2.
Suppose that the maximum $\eta_{\rm max}$ is attained by two points $\xi$ and $\xi'$
which are different modulo ${\rm Im}\, A^T$.
Since $-\log Z$ is an upper convex function,
the function $\eta \cdot \xi - \log Z$ is constant with the value $\eta_{\rm max}$
on the segment $ s \xi + (1-s) \xi'$, $0 \leq s \leq 1$.
Since this function is holomorphic, it is constant on the line defined by $\xi$ and
$\xi'$.
It follows from the proof of 1 that the value of this function is smaller than $M$ when 
$|s \xi + (1-s) \xi'|$ is sufficiently large.
This contradicts that the function is constant on the line.
\qed

{\it Remark} 1.
The existence proof of the maxarg of $f(\xi)$ for $\eta \in {\rm relint}({\rm New}(Z))$
can be easily be extended to a more general model that
$Z=\sum_{u \in S} c_u p^u$,
where $S$ is a finite set in $\N^n$ and $c_u$ is a positive number,
and
$\xi$ is parameterized as $\xi=B \tau$, where $B$ is an $n \times m$ matrix
and $\tau \in {\R}^m$.
\medbreak

{\it Remark} 2. 
We can reduce the proof of our theorem to 
Theorem 2.5 of Haberman \cite{haberman}.
Let us sketch it.
Set $S=S(Z)$.
The theorem of Haberman states that the MLE exists if and only if
$$ S^* = \{ \mu \in {\rm Ker}(A) \,|\,
 (u'-u)\cdot \mu \leq 0 \ \mbox{ for all } u' \in S \}
$$
is $\{ 0 \}$.
When $u$ is a point in ${\rm New}(Z)$, 
the cone $C_u=\{ \mu \,|\, (u'-u)\cdot \mu \leq 0 \ 
 \mbox{for all } u' \in {\rm New}(Z)\}$
is the outer normal cone at $u$.
It is a fundamental result in the theory of polytopes that
the union of these cones is a fan and in particular
$C_u \cap {\rm aff}({\rm New}(Z)) = \{ 0 \}$ 
if and only if $u$ is in the relative interior of the polytope ${\rm New}(Z)$.
Let $u$ be outside the Newton polytope.
Consider the cone $C$ generated by 
$\{u' - u \,|\, u' \in {\rm New}(Z)\}$.
$C$ is strictly contained in the affine hull of ${\rm New}(Z)$.
Therefore, the dual cone of $C$, which is equal to $C_u$,
contains a non-zero vector of the affine hull.
Hence, $S^*$ contains a non-zero vector.
The equivalence of the existence of the MLE and the surjectiveness 
of the moment map can be proved as in the proof of our theorem.
\medbreak

{\it Remark} 3. 
The bijection in the theorem is presented in different forms
in several studies reported in the literature.
We have seen Haberman's result in Remark 2.
Fienberg and Rinaldo \cite{fienberg-rinaldo} give a closely related result for the existence of the 
unconditional MLE, whereas we are concerned with the existence of the conditional MLE.
The recent exciting paper \cite{MSUZ} states that 
``Theorem 2.2 (the image is the convex hull of all sufficient statistics) in this paper is standard in the theory of exponential families
(see \cite[Theorem 3.6]{brown}). 
This paper concerns situations when this bijection has desirable algebraic
properties''.
We characterize the image of the map as the Newton polytope
and will discuss an algorithm for computing the inverse image
by the HGM 
and an asymptotic equivalence of the moment map and IPS.

\bigbreak
The following lemma is used to prove Theorem \ref{th:one-to-one}
and is well known.
We include the proof for the convenience of
readers in the hypergeometric community.
\begin{lemma}  \label{lm:zero}
Let $f(\xi)$ be an upper convex $C^2$ class function.
The existence of the maximum of $f(\xi)$ and the existence of a point $q$
such that ${\rm grad}(f) = 0$ are equivalent.
\end{lemma}

{\it Proof}\/.
We suppose that ${\rm grad}(f)=0$ at $\xi=q$.
We restrict this function to $\xi = q + t \xi^0$,
where $t$ is a scalar variable and $\xi^0$ is any vector.
The function $f(q+t \xi^0)$ is an upper convex function in one variable.
We denote the function by $f(\xi^0; t)$ if no confusion arises.
Since ${\rm grad}(f)=0$ at $\xi=q$, we have $f'(\xi^0; 0)=0$.
We may assume that $f(\xi^0; 0)=0$ without a loss of generality.
Since the upper convexity implies $f''(\xi^0; t) \leq 0$, 
we have
$f'(\xi^0; t) \leq 0$ for $t > 0$ and
$f'(\xi^0; t) \geq 0$ for $t < 0$.
Then, we have $f(\xi^0; t) \leq 0$ in a neighborhood of the origin.
Since a local maximum is the (global) maximum for an upper convex function, 
we have $f(t) \leq 0$ for all $t$.
Since this argument holds for any $\xi^0$,
we have $f(\xi) \leq 0$,
which means that $0$ is the maximum.
The converse is an elementary fact in calculus. \qed

\bigbreak

When the matrix $A$ represents a two-way contingency table,
the image agrees with the open transportation polytope, because in this case
$A$ is totally unimodular.
This follows from the following corollary.
\begin{corollary}
Retain the assumption of the theorem
(the dimension of ${\rm New}(Z)$ is $n-d$).
If $A$ is a totally unimodular matrix, the image of $E[U]$ agrees with
the open transportation polytope.
\end{corollary}

{\it Proof}\/.
Since $A$ is totally unimodular, all the vertices of the transportation polytope
are in $\N^n$.
Then the transportation polytope agrees with
the convex hull of $S(Z)$.
\qed

\begin{example}\rm  \label{ex:2x3}
$$
A = \left(
\begin{array}{cccccc}
1 & 1 & 1 & 0 & 0 & 0\\
0 & 0 & 0 & 1 & 1 & 1 \\
1 & 0 & 0 & 1 & 0 & 0 \\
0 & 1 & 0 & 0 & 1 & 0 \\
\end{array}
\right), \quad
\beta=(b_1,b_2,c_1,c_2).
$$
Then, $u$ satisfying $Au=\beta$ can be regarded as a $2 \times 3$ contingency table
$$ u =
\left(
\begin{array}{ccc}
 u_1 & u_2 & u_3\\
 u_4 & u_5 & u_6 \\
\end{array}
\right)
= 
\left(
\begin{array}{ccc}
 u_{11} & u_{12} & u_{13}\\
 u_{21} & u_{22} & u_{23} \\
\end{array}
\right)
$$
with the row sums $b_1$ and $b_2$
and the column sums $c_1$, $c_2$ and $c_3 = b_1 + b_2 -c_1-c_2$.
When we regard $u$ as a contingency table, we denote
an entry of $u$ with two indices as $u_{ij}$.
Figure \ref{fig:im23} is the image of $(E[U_{11}], E[U_{23}])$
when $(b;c) = (21,7;12,5,11)$.
\begin{figure}[tb]
\setlength\unitlength{4mm}
\begin{picture}(13,8)(4,0)
\put(5,0){\line(1,0){5}}
\put(5,0){$(5,0)$}
\put(10,0){\line(1,1){2}}
\put(10,0){$(10,0)$}
\put(12,2){\line(0,1){5}}
\put(12,2){$(12,2)$}
\put(12,7){\line(-1,-1){7}}
\put(12,7){$(12,7)$}

\end{picture}
\caption{Image of $(E[U_{11}],E[U_{23}])$} \label{fig:im23}
\end{figure}
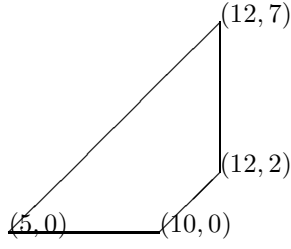
\end{example}

\begin{example}\rm
When $A$ is not totally unimodular, 
the Newton polytope ${\rm New}(Z)$ is not equal to the transportation
polytope in general and the image does not agree with the open transportation polytope.
We consider 
$A = \left(
\begin{array}{ccc}
1 & 1 & 1 \\
0 & 1 & 2 \\
\end{array}
\right)$.
The open transportation polytope is defined by
$$ e_3 > 0,\ \ e_1 = e_3+\beta_1-\beta_2 > 0,\ \ e_2=-2e_3 + \beta_2 >0, $$
which is an open interval. 
When $\beta_2 \leq \beta_1$, we have
$ 0 < e_3 < \beta_2/2$.
Let us specialize to the case $(\beta_1,\beta_2) = (4,3)$.
Then, $Z(\beta;p) = \frac{1}{3!} p_1 p_2^3 + \frac{1}{2!} p_1^2p_2 p_3$.
We have 
\begin{eqnarray*}
p_3 \partial_3 Z/Z
&=& \frac{1/2! p_1^2 p_2 p_3}
 {1/3! p_1 p_2^2 + 1/2! p_1^2 p_2 p_3} \\
&=& \frac{1}{1+\frac{2!}{3!} \frac{p_2^2}{p_1 p_3}}.
\end{eqnarray*}
Hence, the expectation of $U_3$ agrees with the interval $(0,1)$
which is contained in $(0,\beta_2/2)=(0,3/2)$ and does not agree with it.
\end{example}

\section{Classical Hypergeometric Polynomials} \label{sec:classical}

We consider 
the $(2+m) \times 2m$ matrix 
$$A = \left(
\begin{array}{cccccccccccc}
1 & 1&1& \cdots &1 & 1 & 0 & 0 &0& \cdots &0 & 0 \\
0 & 0&0& \cdots &0 & 0 & 1 & 1 &1& \cdots &1 & 1 \\
1 & 0&0& \cdots &0 & 0 & 1 & 0 &0& \cdots &0 & 0 \\
0 & 1&0& \cdots &0 & 0 & 0 & 1 &0& \cdots &0 & 0 \\ 
 &&&\cdots&&&&&&\cdots&& \\
0 & 0&0& \cdots &1 & 0 & 0 & 0 &0& \cdots &1 & 0 \\ 
\end{array}
\right)
$$ 
representing a $2 \times (m+1)$ contingency table.
When the first column sum is greater than or equal to the first row sum,
the normalizing constant is expressed in terms of the Lauricella function $F_D$
(see, e.g., \cite{goto}),
which is defined by
$$ F_D(a,b,c;z_1, \ldots, z_m) =
\sum_{k \in \N^m} 
 \frac{(a)_{|k|} (b_1)_{k_1} \cdots (b_m)_{k_m}}
      {(c)_{|k|} (1)_{k_1} \cdots (1)_{k_m}} z^k.
$$
Here, $a$, $b=(b_1, \ldots, b_m)$, and $c$ are numbers, 
$|k| = k_1 + \cdots +k_m$, $(\alpha)_i = \alpha (\alpha+1) \cdots (\alpha+i-1)$
is the Pochhammer symbol, and
$z^k = \prod_{i=1}^m z_i^{k_i}$.
For non-positive integers $a$ and $b_1, \ldots, b_m$ and a positive integer $c$,
set 
$$\mu = \left(
\begin{array}{ccccc}
-a & 0  & 0  & \cdots & 0    \\
c-1&-b_1&-b_2& \cdots & -b_m \\
\end{array}
\right).
$$
Then, the normalizing constant $Z(\beta;p)$ for
$\beta=(-a, c-1-\sum b_i,c-1-a, -b_1, \ldots, -b_{m-1})$
is equal to
$$ p^\mu F_D(a,b,c; z), \quad z_i = \frac{p_{0i}p_{10}}{p_{00} p_{1i}},\ \ i=1, \ldots, m. $$
Here, 
$p = \left(
\begin{array}{cccc}
p_{00} & p_{01} & \cdots &p_{0m} \\
p_{10} & p_{11} & \cdots &p_{1m} \\
\end{array}
\right)
$ and is regarded as a vector of length $2(m+1)$.
Our Proposition \ref{prop:invariance} and Theorem \ref{th:one-to-one} yield the following.
\begin{theorem}  \label{th:fdimage}
Assume that $a$ and $b_i$, $i=1, \ldots, m$,
are non-positive integers and that $c$ is a positive integer.
The image of the map
$$ \left( z_1 \frac{\partial F_D(a, b,c; z)}{\partial z_1}/F_D(a,b,c;z),
             \ldots,
             z_m \frac{\partial F_D(a, b,c; z)}{\partial z_m}/F_D(a,b,c;z)
     \right)
$$
for $z \in {\R}_{>0}^m$
is the relative interior of the polytope defined by
$$
 \{ \eta \in {\R}^m \,|\, \eta_1 + \cdots + \eta_m \leq -a,
    \eta_i \leq -b_i \mbox{ and } \eta_i \geq 0 \mbox{ for $i=1, \ldots, m$}\},
$$
and the map is one-to-one.
\end{theorem}

{\it Proof}\/.
When the Newton polytope has the dimension $m$, the conclusion follows from 
Proposition \ref{prop:invariance} and Theorem \ref{th:one-to-one}
and the expression of $Z$ in terms of $F_D$.
When the dimension of the Newton polytope is smaller than $m$,
$F_D$ is expressed in terms of $F_D$ of a smaller set of independent variables 
or a constant
and at least one of $a$, $b_1, \ldots, b_m$ is $0$.
We can prove the conclusion by case by case checks.
\qed

From the viewpoint hypergeometric functions,
we can understand that the expectation is expressed as a quotient of
the hypergeometric polynomial $Z(\beta;p)$ and its derivative.
It is natural to consider a moment map as representing a basis of the twisted cohomology group
to a projective space.
This moment map is called the co-Schwartz map in recent studies
of the hyperbolic Schwartz map \cite{SY}.

Once the normalizing constant $Z$ is expressed in terms of the Lauricella function $F_D$,
we can utilize several formulas in the study of hypergeometric functions
for efficient numerical evaluation of the normalizing constant and its derivatives
by the HGM.
For details, see \cite{goto}, \cite{ott}, and the Risa/Asir package
{\tt tk\_fd.rr} \cite{asir}.
Our examples in this paper are calculated with this package.

\section{MLE Algorithms and Information Geometry}  \label{sec:mle}

For given data $u \in \N^n$ and configuration matrix $A$,
we are interested in numerically solving the conditional MLE problem
$$ {\rm maxarg}_\xi\,\frac{\exp(\xi)^u}{u! Z(\beta; \exp(\xi))}, \quad \beta=Au. $$
Taking the logarithm of the likelihood function, we have
\begin{equation} \label{eq:loglikely}
 u \cdot \xi - \log u! - \log Z. 
\end{equation}
The maximization of this function is equivalent to the maximization of (\ref{eq:loglike})
with $\eta = u$.
In \cite{ott}, we solve this MLE problem with the 
HGM with respect to discrete parameters $\beta$
and also show that the BFGS algorithm 
should be used for the optimization part.
We explain a method for solving the MLE problem 
by using the framework of information geometry.
The maximization problem can be transformed into a problem of solving 
a system of algebraic equations.
In fact,
differentiating (\ref{eq:loglikely}) by $\xi_i$, we obtain
$ u_i = \frac{\partial Z}{\partial \xi_i} / Z $.
We use the variable $p_i = \exp(\xi_i)$ to present the method.
Recall that $A$ is a $d \times n$ matrix which defines the 
$A$-hypergeometric distribution.
In information geometry, the space of $\xi$ and the space of $\eta$ are dual.
We give a one-to-one correspondence in Theorem \ref{th:one-to-one}.
We assume, for simplicity of presentation, 
$$ {\R}_{>0} \ni y=(p_{d+1}, \ldots, p_n) \ \mapsto\  
 E(y)=(E[U_{d+1}], \ldots, E[U_n]) \in {\rm relint}({\rm New}(Z)) \cap {\R}^{n-d}
$$
is a one-to-one correspondence when $p_1, \ldots, p_{d}$ are fixed.
In other words, when $p_i=\exp(\xi_i)$
and $\xi_1, \ldots, \xi_d$ are fixed,
$\xi_{d+1}, \ldots, \xi_n$ are complete representatives of $\R^n/{\rm Im}\, A^T$.
We are interested in an algorithm for finding $y^*$ for a given $\eta^*$
such that $E(y^*)=\eta^*$.
By a recipe in information geometry  or by Newton's method, 
the inverse point $y^*$
of $\eta^*$ can be obtained by the iteration
\begin{equation} \label{eq:iteration}
 \mbox{ new $y$} = y +
  \varepsilon {\dot E}(y)^{-1} (\eta^*-\eta),
\end{equation}
where 
${\dot E}(y) = \left( \frac{\partial E[U_{d+i}]}
                           {\partial p_{d+j}} \right)_{i,j = 1, \ldots, n-d}
$,
$\eta = E(y)$, and $\varepsilon$ is a (sufficiently small) number.
This iteration formula can be obtained as follows.
In the information geometry algorithm,
the inverse image of
the movement from $\eta$ to $\eta^*$ along a straight line will give
a good movement  in the $y$-space.
We consider the first-order approximation of $E$ as follows.
$$ E(y+h) \sim E(y) + {\dot E}(y) h, \quad
  {\dot E}(y) = \left(
\begin{array}{cccc}
 \frac{\partial E[U_{d+1}]}{\partial p_{d+1}}& 
 \frac{\partial E[U_{d+1}]}{\partial p_{d+2}}& 
 \cdots &
 \frac{\partial E[U_{d+1}]}{\partial p_{n}} \\ 
 \frac{\partial E[U_{d+2}]}{\partial p_{d+1}}& 
 \frac{\partial E[U_{d+2}]}{\partial p_{d+2}}& 
 \cdots &
 \frac{\partial E[U_{d+2}]}{\partial p_{n}} \\ 
 &&\cdots& \\
 &&\cdots& \\
 \frac{\partial E[U_{n}]}{\partial p_{d+1}}& 
 \frac{\partial E[U_{n}]}{\partial p_{d+2}}& 
 \cdots &
 \frac{\partial E[U_{n}]}{\partial p_{n}} \\ 
\end{array}
\right) 
$$
If $E(y+h)=\eta^*$ and $E(y)=\eta$, 
then $h$ is approximately equal to
${\dot E}(y)^{-1} (\eta^*-\eta)$.
Thus, we may expect that $h$ is a good direction for updating $y$ to a new $y$.

The gradient matrix ${\dot E}(y)$ can be evaluated by the Pfaffian system
for the HGM
\cite{N3OST2}, \cite{ot2}.
Let us briefly summarize the evaluation method.
We regard the $\beta_j$'s as indeterminates in the following discussion.
Let $s_1=1, s_2, \ldots, s_r \in D$ be the standard monomials of the 
$A$-hypergeometric system for $Z$.
Here $D$ is the ring of differential operators
${\bf Q}(\beta_1, \ldots, \beta_d)\langle p_1, \ldots, p_n,
\pd{p_1}, \ldots, \pd{p_n} \rangle$,
where $\pd{p_i} = \partial/\partial p_i$.
Let $F=(Z, s_2 \bullet Z, \ldots, s_r \bullet Z)^T$.
Then, the vector valued function $F$ satisfies the Pfaffian system
$$ \frac{\partial F}{\partial p_i} = P_i(\beta,p) F, \quad
 i =1, \ldots, n,
$$
where $P_i$ is an $r \times r$ matrix with rational function entries
with respect to $\beta$ and $p$.
Differentiating both sides of the Pfaffian system by $p_j$,
we have
\begin{equation} \label{eq:hessianPf}
 \frac{\partial^2 F}{\partial p_i \partial p_j} =
  \frac{\partial P_i(\beta,p)}{\partial p_j} F +
  P_i \frac{\partial F}{\partial p_j} = 
  \left( \frac{\partial P_i(\beta,p)}{\partial p_j} +
  P_i P_j \right) F.
\end{equation}
Therefore, the numerical value of the left-hand side can be 
evaluated from the numerical value of $F$.
The $k$-th entry of 
$\frac{\partial F}{\partial p_i}$ is
$\pd{i} s_k \bullet Z$ and
the $k$-th entry of 
$\frac{\partial^2 F}{\partial p_i \partial p_j}$
is equal to
$\pd{p_i}\pd{p_j}s_k \bullet Z$.
Since $E[U_i](p) = p_i \pd{p_i} \bullet \log Z$, 
$\pd{j}\bullet E[U_i](p) 
= \pd{p_j} \bullet (p_i (\pd{p_i}\bullet Z)/Z)$.
\comment{
We suppose that $\pd{p_{d+1}}, \ldots, \pd{p_{n}}$ are contained
in the set $s_1, \ldots, s_r$
for simplicity.
Under this assumption, 
} Thus, the numerical value of $\pd{j}\bullet E[U_i](p)$ can be obtained 
from the numerical value of $F$, which can be evaluated by
the discrete HGM \cite{ot2}.
Let us next discuss the convergence of our method.
\begin{proposition}  \label{prop:negative_definite}
Under the assumption of Theorem \ref{th:one-to-one}
and the assumption that
the $(\xi_{d+1}, \ldots, \xi_n)$'s are complete representatives of
$\R^n / {\rm Im}\, A^T$
for fixed $\xi_1, \ldots, \xi_d$,
the matrix ${\dot E}(y)$ is a negative definite matrix for any $y \in \R^{n-d}_{>0}$.
\end{proposition}

{\it Proof}\/.
Set $H(t) = \sum_{i=1}^t \sum_{k=M}^{M+N} c_{ik} \exp(k \alpha_i t)$,
where
$c_{ik} \geq 0$ and $\alpha_i \in \R$ are linearly independent over $\Z$. 
We assume that $H(t)$ has at least two non-zero terms.
From the proof of Lemma \ref{lm:exppoly},  we have
$\frac{d^2}{dt^2} \log H(t) > 0$.

Fixing real numbers $\delta_{d+1}, \ldots, \delta_n$
and
$\gamma_{d+1}, \ldots, \gamma_n$,
we restrict $-\log Z$ to
$q(t)=(\exp(\xi_1), \ldots, \exp(\xi_d),
  \exp(\delta_{d+1} t + \gamma_{d+1}), \ldots,  
  \exp(\delta_{n} t + \gamma_{n}))$.
Set 
$F(t) = -\log Z(\beta; q(t))$.
Since the dimension of the Newton polytope of $Z$ is $n-d$,
the Newton polytope of $Z(\beta;\exp(\xi_1), \ldots, \exp(\xi_d),
p_{d+1}, \ldots, p_n)$ as the polynomial in $p_{d+1}, \ldots, p_n$
is $n-d$-dimensional from the assumption.
Therefore, 
the restricted $Z$ has at least two non-zero terms as a polynomial in
$\exp(n\cdot \delta t)$, $n \in \Z^n$,
where $\delta=(0, \ldots, 0,\delta_{d+1}, \ldots, \delta_n)$.
From the observation at the beginning of the proof, we have
$\frac{d^2}{dt^2} F(t) < 0$.

Let $G(\xi)$ be the Hessian of $-\log Z(\beta;p(\xi))$
with respect to $\xi_{d+1}, \ldots, \xi_n$.
We have
$G=
\left( - \frac{\partial}{\partial \xi_{d+i}}
\left( \frac{\partial Z/\partial \xi_{d+j}}{Z} \right)\right)
$,
$F'(t) = \frac{-1}{Z} \sum_{j=1}^{n-d} \frac{\partial Z}{\partial \xi_{d+j}}
 \delta_{d+j}$,
and
$F''(t) = \sum_{i=1}^{n-d} \delta_{d+i} \frac{\partial}{\partial \xi_{d+i}}
\left( - \sum_{j=1}^{n-d} \frac{\partial Z/\partial \xi_{d+j}}{Z}
 \delta_{d+j} \right)
= \sum_{i,j=1}^{n-d} G_{ij} \delta_{d+i} \delta_{d+j}
$.
Assume that $G$ has a positive or zero eigenvalue $\lambda \geq 0$
at $(\xi_{d+1}, \ldots, \xi_n) = \gamma$.
Letting $\delta$ be an eigenvector for $\lambda$,
we then have $\delta^T G \delta = \lambda |\delta|^2 \geq 0$.
This contradicts that $F'' < 0$.
Therefore, the matrix $G$ is negative definite.
Since ${\rm diag}(p_{d+1}, \ldots, p_n) {\dot E}(y) = G$,
the matrix ${\dot E}(y)$ is negative definite.
\qed

It is well known that if ${\dot E}(y)$ is negative definite at $y=y^*$,
the iteration (\ref{eq:iteration}) converges when the starting $y$ is 
sufficiently close to $y^*$
(see, e.g., \cite[Th 3.5]{nocedall}).
Then, the remaining task we need to do is to find $y$ which is sufficiently
close to $y^*$.
This problem will be discussed in the next two sections.
Let us briefly summarize it.
For observed data $u$,
we take $u/|u|$ ($|u|=u_1+\cdots+u_n$) 
as the initial value of the iteration (\ref{eq:iteration}). 
This choice is expected to work well,
because when $p=u/|u|$, the approximate expectation evaluated
by the IPS output is close to $u$, as we will see in sections \ref{sec:IPS}
and \ref{sec:equivalence}.
The following example illustrates the effectiveness of our method.

\begin{example} \rm

Let $A$ be the matrix
$$\left(
\begin{array}{ccccccc}
 0 & 0 & 0 & 1 & 1 & 1 & 1 \\
 1 & 0 & 0 & 1 & 0 & 1 & 0 \\
 0 & 1 & 1 & 0 & 1 & 0 & 1 \\
 1 & 1 & 0 & 1 & 1 & 0 & 0 
\end{array}
\right)
$$
and consider the discrete $A$-hypergeometric distribution
defined by this $A$.
The model defined by $A$ can be regarded as 
a $2 \times 2 \times 2$ contingency table with one structural $0$
with fixed one-dimensional marginal sums.
In other words, the model represents the table
$$
\begin{array}{c|cc}
 && \\ \hline
 & p_1 & p_2 \\
 & 0   & p_3 \\
\end{array}
\quad
\begin{array}{c|cc}
&& \\ \hline
& p_4 & p_5 \\
& p_6 & p_7 \\
\end{array},
$$
with fixed marginal sums of ``planes'' of the cube
$p_4+p_5+p_6+p_7$,
$p_1+p_4+p_5$,
$p_2+p_3+p_5+p_7$,
$p_1+p_2+p_4+p_5$.
Assume that we observe the data $\eta^*$,
$$
\begin{array}{c|cc}
&& \\ \hline
& 19 & 132 \\
& 0  & 9 \\
\end{array}
\quad
\begin{array}{c|cc}
&& \\ \hline
& 11 & 52 \\
&  6 & 97 \\
\end{array}
$$
The total number of incidences is 
$19+132+9+11+52+6+97=326$.
We want to find an approximate value of 
$p$ such that the vector $(E[U_i])$ agrees with the observed data
$\eta^*$.
As the first approximation of $p$, 
we take $P_0 = (19,132,9,11,52,6,97)/326$.
The rank $r$ of the corresponding $A$-hypergeometric system is $5$ and
we can take the set of standard monomials 
$\{s_i\} = \{1, \pd{p_5},\pd{p_6}, \pd{p_7}, \pd{p_7}^2\}$.
Note that the expectation polytope ${\rm New}(Z)$ is $3$-dimensional and 
there exists a one-to-one correspondence between 
$(p_5,p_6,p_7)$ and $\eta=(E[U_5],E[U_6],E[U_7])$
when $p_1, p_2, p_3, p_4$ are fixed.
Our HGM software evaluates the expectation at $P_0$ with rational 
arithmetic \cite{ot2}. 
The approximate value of $\eta$ is
$( 51.9194, 5.99193, 97.0891 )$.
This value is close to $\eta^*=(52,6,97)$
and the error is bounded by $0.09$.
We refine the value $P_0$ by the gradient of the expectation ${\dot E}(y)$
evaluated by the HGM and the derivative of the Pfaffian system (\ref{eq:hessianPf}).
The value of $({\dot E})^{-1} (\eta^*-\eta)$ is approximately equal to
$h=(0.000256154, -0.000152585, -0.00310983)$.
As explained in (\ref{eq:iteration}),
we update $P_0$ and define the new $P_0$, which is denoted by $P_1$,
as $P_0+(0,0,0,0,h_1,h_2,h_3)$.
These steps are performed in 11.1s by Risa/Asir
on a machine with an Intel Xeon
CPU (2.70 GHz) and 256 G of memory.

We again apply the HGM and evaluate the expectation and its gradient.
The approximate value of $\eta$ for $P_1$ is
$\eta=( 52.0006, 6.00006, 96.9993)$
and the error is bounded by $0.0007$ ($\ll 0.09$).
The second step takes 236 s, because we have big denominators and numerators
in the rational arithmetic calculation. 
We have obtained a very good approximation of the MLE of $p$
by only two iterations.
\end{example}

\section{Generalized Odds Ratio and IPS}
\label{sec:IPS}

As a preliminary to the asymptotic analysis of the next section, we discuss the generalized odds ratio, log-affine model, and IPS.

Let $\A$ 
be the Gale transform of $A$ as defined in section \ref{sec:expectation}.
From the construction, 
$\A$ has the properties that $\A$ is a full-rank matrix such that $A \A^T =0$,
and that any $u\in\N^n$ such that $Au=\beta$ can be written as $u=u_0+\A^T w$, where $u_0\in\N^n$ is a fixed point such that $Au_0=\beta$ and $w=(w_1,\ldots,w_{n-d})^{T}\in\Z^{n-d}$.
%
The generalized odds ratio has been defined as
$\exp(\lambda_i) = p^{\a_i} = p_1^{\a_{i1}}\cdots p_n^{\a_{in}}$, $i=1,\ldots,n-d$,
and we define the generalized log odds ratio $\lambda$ as
\[
 \lambda=(\lambda_1,\ldots,\lambda_{n-d})^T=\A\xi.
\]
It is easy to see that the parameter $\exp(\lambda)$ or $\lambda$ is one-to-one to the set of probability distributions (\ref{eq:a-distribution}).

For the $A$ in Example \ref{ex:2F1a}, we can choose $\A=(1,-1,-1,1)$, and 
$\exp(\lambda) = p^{(1,-1,-1,1)} = p_1 p_2^{-1} p_3^{-1} p_4$ or
$\lambda = \log (p_1 p_4/ p_2 p_3)$.
This is nothing but the classical (log) odds ratio of a $2\times 2$ table.

For the $A$ in Example \ref{ex:2x3}, we can choose
\[
 \A = \begin{pmatrix}
 1 & -1 & 0 & -1 & 1 & 0 \\
 1 & 0 & -1 & -1 & 0 & 1
 \end{pmatrix}
\ \ \mbox{and}\ \ %
 \lambda = \left(\log\frac{p_1 p_5}{p_2 p_4},\log\frac{p_1 p_6}{p_3 p_4}\right)^T.
\]

We consider the Poisson distribution according to
$ \frac{p^u}{u!}\exp(-{\bf 1} \cdot p)$
with the affine structure
\begin{equation}
\label{eq:log-affine}
 \log p \in \A^T(\A\A^T)^{-1}\lambda + {\rm Im}\,A^T.
\end{equation}
Here $\lambda=\log p^{\A}=(\lambda_1,\ldots,\lambda_{n-d})^T$ is assumed to be fixed.
This statistical model is called the log-affine model \cite{lauritzen}.
Let $U$, which takes values in $\N^n$, 
be distributed according to this log-affine model
parameterized  by the column vector $\theta \in \R^d$ as
$$
p(\theta) = \exp\left(
  \A^T(\A\A^T)^{-1}\lambda + A^T \theta \right),
$$
where $\exp(v) = (\exp(v_1), \ldots, \exp(v_n))^T$ for $v \in \R^n$.
Note that it is the unconditional model (\ref{eq:unconditional})
with a different parameterization $p(\theta)$ above 
from that by $\xi$ in the previous sections.
In the log-affine model, we do not impose the condition $AU=\beta$.
Let $\theta^*$ be the MLE of $\theta$ and let $m=p(\theta^*)$ if they exist.
The log-affine model is an exponential family with sufficient statistics $T=AU$.
The convex hull of the support of $T$ is $K=\mathrm{conv}(\N A)=\R_{\ge 0} A$.
According to the general theory of the exponential family, if $T\in\mathrm{int}\,K$, then the MLE exists uniquely (\cite[Theorem 9.13]{barndorff-nielsen},\cite[Theorem 3.6]{brown}).
In this case,
$m=(m_1,\ldots,m_n)^T = p(\theta^*)$ is the unique solution of
\[
\left\{
\begin{split}
 & A U = A m, \\
 & \lambda 
 = \A \log m.
\end{split}\right.
\]
In particular, $m_i>0$.
Note that $m$ is a function of $A U$ and $\lambda$.
IPS is a numerical procedure for obtaining $m$ when $A U$ and $\lambda$ are given.
Although IPS was originally invented for contingency tables and hierarchical models (\cite{sinkhorn-knopp,darroch-ratcliff, lauritzen}),
this procedure can be extended to the log-affine model as follows.
\begin{theorem}
\label{th:ips}
The IPS for the log-affine model (\ref{eq:log-affine}):
\begin{enumerate}
\item
 Set $m^{(0)} := \exp(\A^T(\A\A^T)^{-1}\lambda)$ as an initial value.
\item
 For $t=0,1,2,\ldots$, let
\[
 m_i^{(t+1)} := m_i^{(t)} \exp(\mu^{(t)}a_{ji}), \ \ j=(t\ \mathrm{mod}\ d)+1, \ \ i=1,\ldots,n,
\]
where $\mu^{(t)}$ is the unique solution of
\[
 \sum_{i=1}^n a_{ji} m_i^{(t)} \exp(\mu^{(t)}a_{ji}) = (AU)_j,
\]
and $()_j$ is the $j$-th element of a vector.

\item
The limit $m=\lim_{t\to\infty} m^{(t)}$ is the desired output of IPS.
\end{enumerate}
\end{theorem}

{\it Proof\/}.
IPS is interpreted as a method for solving the dual problem of maximizing likelihood (\cite{csiszar,dykstra-lemke}).
For $p,q\in\R_{>0}^n$, define the $I$-divergence
\[
 I(p\Vert q) = \left( -\log\frac{q}{p} +\frac{q}{p} -{\bf 1} \right) \cdot p,
\]
where $q/p=(q_1/p_1,\ldots,q_n/p_n)$.
Then, the MLE $m$ is the minimizer of the minimizing problem
\[
 \mbox{\tt Minimize} \ \ I(U\Vert q)\ \ \mbox{\tt subject to}\ \ \log q\in\xi_0+L,
\]
where $\xi_0=\A^T(\A\A^T)^{-1}\lambda$, $L=\mathrm{Im}A^T$.
The variable $q$ is parameterized as $\exp(\xi_0 + A^T \theta)$.
This is a convex problem and its dual problem gives the same answer.
The dual problem is formalized as
\[
 \mbox{\tt Minimize} \ \ I(p\Vert p_0)\ \ \mbox{\tt subject to}\ \ p\in U+L^\perp,
\]
where $p_0=\exp(\xi_0)$.
Here $L^\perp=\mathrm{Ker}A = \bigcap_{i=1}^d M_i$, 
and the linear space
$M_i$ is the orthogonal complement of the $i$-th row vector of the matrix $A$.
Starting from $m^{(0)}=p_0$, IPS is the procedure for generating a sequence $m^{(t)}$, $t=0,1,\ldots$, by
\[
 m^{(t+1)}:= \mathrm{argmin}_{m\in U+M_j} I(m\Vert m^{(t)}), \ \ j=(t\ \mathrm{mod}\ d)+1.
\]
Noting that $m\in U+M_j$ $\Leftrightarrow$ $\sum_{i=1}^n a_{ji}m_i = (AU)_j$, and that
$I(m\Vert m^{(t)})=\sum_{i=1}^n m_i(\log (m_i/m^{(t)}_i) -1) + \mathrm{const}$, we define the Lagrangian
\[
 L = \sum_{i=1}^n m_i\left( \log\frac{m_i}{m^{(t)}_i} -1 \right) - \mu \left( \sum_{i=1}^n a_{ji}m_i - (AU)_j \right),
\]
and
$0=\partial L/\partial m_i = \log(m_i/m^{(t)}_i)-\mu a_{ji}$ yields $m_i=m_i^{(t)}\exp(\mu a_{ji})$.
$\mu$ is determined by $\sum_{i=1}^n a_{ji}m_i = (AU)_j$.
\qed

\section{Asymptotic Equivalence of the Moment Map and IPS}
\label{sec:equivalence}

The $m$ by IPS will be a good approximation of the expectation
of $u$ of the $A$-hypergeometric distribution.
This fact is well known for some contingency tables as illustrated in,
e.g., books \cite{hirotsu}, \cite{plackett}.
We will show this fact for any $A$-hypergeometric distribution by
extending and validating Plackett's
\cite[pp. 41 ($2\times 2$ table), pp. 65--66 ($r\times s$ table)]{plackett}
and Hirotsu's \cite{hirotsu} heuristic idea, borrowing techniques from the local central limit theorem
\cite[Section I.6, pp.56)]{shiryaev}.

We define a series of probability distributions
\[
 P_\kappa(u,\xi) = \frac{\exp(u \cdot \xi)}{u! Z_\kappa(\xi)}, \quad u\in S_\kappa, \quad \kappa=1,2,\ldots,
\] 
where
\[
 S_\kappa = \{ u\in\N^n \,|\, A u = \kappa\beta \}, \quad
 Z_\kappa(\xi) = \sum_{u\in S_\kappa} \frac{\exp(u\cdot\xi)}{u!},
\]
and consider its limiting behavior when $\kappa\to\infty$.

Let $m=m(\lambda)$ be the unique solution of the following IPS: 
\begin{equation}
\label{eq:ips}
\left\{
\begin{split}
 & \beta = A m, \\
 & \lambda 
 = \A \log m.
\end{split}\right.
\end{equation}
As explained in the previous section, if $\beta\in\mathrm{int}(\R_{\ge 0}A)$, $m>0$ is determined uniquely.
We will establish an asymptotic approximation of the $A$-hypergeometric distribution to a Gaussian density function.
\begin{theorem}
\label{th:asymptotic}
Suppose that $\beta\in\N A\cap\mathrm{int}(\R_{\ge 0}A)$.
\begin{equation*}
\sup_{\forall i\,|u_i-\kappa m_i|<\varphi(\kappa)} \left\vert\frac{P_\kappa(u,\xi)}{\widehat P_\kappa(u,\xi)}-1\right\vert \to 0
\ \ (\kappa\to\infty),
\end{equation*}
where $\varphi(\kappa)$ is a positive function satisfying 
$\varphi(\kappa)=o(\kappa^{2/3})$, 
$\kappa/\varphi(\kappa)^2 = o(1)$,
and
\[
 \widehat P_\kappa(u,\xi) =
 \frac{\det(\A M^{-1}\A^T)^{1/2}}{(2\pi\kappa)^{(n-d)/2}
 }
 \exp\left(-\sum_{i=1}^n \frac{(u_i-\kappa m_i)^2}{2\kappa m_i}\right)
\]
with $M=\mathrm{diag}(m_i)$.
\end{theorem}
We can choose the function $\varphi(\kappa)=\kappa^{7/12}$, for example.
Note that $\widehat P_\kappa(u,\xi)$ depends on $\xi$ through $\lambda=\A\xi$ as expected, since
$m_i=m_i(\lambda)$ depends on $\lambda$.

\medskip
{\it Proof of Theorem \ref{th:asymptotic}\/}.
Write $v=u-\kappa m$.
We first examine the density ratio
\begin{equation*}
 \log\frac{P_\kappa(u,\xi)}{P_\kappa(\kappa m,\xi)}
= (u-\kappa m)\cdot\xi - \log\frac{u!}{(\kappa m)!} \\
= v\cdot\xi - \log\frac{(\kappa m+v)!}{(\kappa m)!}.
\end{equation*}
Here, $(\kappa m)!$ means $\prod_{i=1}^n \Gamma(\kappa m_i + 1)$.
By Stirling's formula
\begin{equation}
\label{eq:stirling}
 \log u! = u(\log u -1) + \frac{1}{2}\log(2\pi u) + R(u),
\end{equation}
where $R(u)=o(1)$ as $u\to\infty$, we have
\begin{align}
 \log\frac{(\kappa m_i+v_i)!}{(\kappa m_i)!}
=& v_i \log(\kappa m_i) + (\kappa m_i) H\left(\frac{v_i}{\kappa m_i}\right) + \frac{1}{2}\log\left(1+\frac{v_i}{\kappa m_i}\right) \nonumber \\
 & + R(\kappa m_i+v_i) - R(\kappa m_i)
\label{eq:ratio_i}
\end{align}
for $\kappa m_i + v_i \geq 1$,
where
\[
 H(v) = (1+v)\log(1+v)-v.
\]

Because of the assumption that $A$ is a configuration matrix, which means
$(1,\ldots,1)^T_{n\times 1} \in {\rm Im}\, A^T$,
and $A v=Au-\kappa Am=0$,
we have
\begin{equation}
\label{eq:sum_v}
 \sum_{i=1}^n v_i=0.
\end{equation}
Moreover, noting that $A v=0 \Leftrightarrow v=\A^T w$, $\exists w=(w_1,\ldots,w_d)^T$, we have
$v^T\log m = w^T\A\log m = w^T\A\log p = v^T\log p = v^T\xi$,
and hence
\begin{equation}
\label{eq:sum_vlogm}
 \sum_{i=1}^n v_i \log m_i = \sum_{i=1}^n v_i \xi_i.
\end{equation}
$H(v)$ has Taylor's expansion
\begin{equation}
\label{eq:taylor}
 H(v) = \frac{1}{2}v^2 - \frac{v^3}{6(1+\theta v)^2}, \ \ 0 < \theta < 1.
\end{equation}
Substituting (\ref{eq:sum_v}), (\ref{eq:sum_vlogm}), and (\ref{eq:taylor}) into (\ref{eq:ratio_i}), and by summing with respect to $i$, we have
\begin{align*}
 \log\frac{P_\kappa(u,\xi)}{P_\kappa(\kappa m,\xi)}
=& -\sum_{i=1}^n\frac{v_i^2}{2 \kappa m_i} \\
 & + \sum_{i=1}^n \left\{ \frac{v_i^3/(\kappa m_i)^2}{6(1+\theta_i \frac{v_i}{\kappa m_i})^2} -\frac{1}{2}\log\left(1+\frac{v_i}{\kappa m_i}\right) -R(\kappa m_i+v_i) +R(\kappa m_i) \right\}
\end{align*}
with $ 0 <\theta_i < 1$.
As $\kappa\to\infty$, the remainder term is $o(1)$ when $v_i=o(\kappa^{2/3})$.
Hence, we have
\begin{equation} \label{lm:asymptotic1}
\sup_{\forall i,\,|u_i-\kappa m_i|<\varphi(\kappa)} \left\vert\frac{P_\kappa(u,\xi)}{P_\kappa(\kappa m,\xi)\exp\left(-\sum_{i=1}^n \frac{(u_i-\kappa m_i)^2}{2 \kappa m_i}\right)}-1\right\vert \to 0
\ \ (\kappa\to\infty),
\end{equation}
where $\varphi(\kappa)=o(\kappa^{2/3})$. 

To complete the asymptotic evaluation of $P_\kappa(u,\xi)$,
we need to evaluate $P_\kappa(\kappa m,\xi)$ as $\kappa\to\infty$.
Note first that
\begin{equation}
\label{eq:summation}
 P_\kappa(\kappa m,\xi)^{-1} = \sum_{u\in S_\kappa}\frac{e^{(u-\kappa m)\cdot\xi} (\kappa m)!}{u!}.
\end{equation}
By letting $u=\kappa m+v$ again, using Stirling's formula, and conducting the same calculations as before, we see that
\begin{equation}
\label{eq:asymptotic2}
\sup_{\forall i,\,|u_i-\kappa m_i|<\varphi(\kappa)} \left\vert\frac
{e^{(u-\kappa m)\cdot\xi} (\kappa m)!/u!}{\exp\left(-\sum_{i=1}^n \frac{(u_i-\kappa m_i)^2}{2 \kappa m_i}\right)}
-1\right\vert \to 0
\ \ (\kappa\to\infty),
\end{equation}
where $\varphi(\kappa)=o(\kappa^{2/3})$.

Let $u_0\in S_1$ be fixed.
The set $S_\kappa$ 
can be written as
\[
 S_\kappa = \{ u\ge 0 \,|\, u= \kappa u_0 + \A^T w,\ w\in \Z^{n-d} \}.
\]
From this observation as well as the uniform approximation given by (\ref{lm:asymptotic1}), the summation over $u\in S_\kappa$ such that $|u_i-\kappa m_i|<\varphi(\kappa)$ can be approximated by the Riemann integral
\begin{align}
&
 \sum_{\forall i,|u_i-\kappa m_i|<\varphi(\kappa)}
 \frac{e^{(u-\kappa m)\cdot\xi} (\kappa m)!}{u!} \nonumber \\
&= (1+o(1))
 \int_{\forall i,|(\A^T w)_i|<\varphi(\kappa)} \exp\left(-\sum_{i=1}^n \frac{(\A^T w)_i^2}{2\kappa m_i}\right) dw
 + O(\kappa^{(n-d-1)/2})
\label{eq:riemann}
\end{align}
as $\kappa\to\infty$.
Moreover, if $\kappa/\varphi(\kappa)^2=o(1)$,
(\ref{eq:riemann}) is asymptotically equivalent to
\begin{equation}
\label{eq:normalizing}
\int_{\R^{n-d}}
 \exp\left(-\sum_{i=1}^n \frac{(\A^T w)_i^2}{2\kappa m_i}\right) dw
 = (2\pi\kappa)^{(n-d)/2} 
\frac{1}{\det(\A M^{-1}\A^T)^{1/2}}
\end{equation}
with $M=\mathrm{diag}(m_i)$,
because, by making the change of variable $w'=w/\sqrt{\kappa}$, the range of integration 
$\{w'\in\R^{n-d} \,|\, \forall i,|(\A^T w')_i|<\varphi(\kappa)/\sqrt{\kappa}\}$ goes to the whole space $\R^{n-d}$.

Next we will see that in the summation (\ref{eq:summation}), the contribution of the outside of $|u_i-\kappa m_i|<\varphi(\kappa)$ is 
negligible. 
Recall that in the Stirling's formula (\ref{eq:stirling}), the upper and lower bounds for the remainder is available (\cite{robbins}):
\[
  \frac{1}{12u+1}<R(u)<\frac{1}{12u}\quad \mbox{for $u\ge 1$}.
\]

Suppose first the case $u_i=\kappa m_i+v_i\ge 1$.
From (\ref{eq:ratio_i}) and the inequality $H(v)\ge H(|v|)$,
$\log(\kappa m_i+v_i)!/(\kappa m_i)!$ is bounded below by
\[
 v_i \log(\kappa m_i) - \frac{1}{2}\log(\kappa m_i) + (\kappa m_i) H\left(\frac{|v_i|}{\kappa m_i}\right) - \frac{2}{12}.
\]
Moreover, when we take a suitable $\kappa_0$, 
the third term is bounded below 
for all $\kappa\ge \kappa_0$
as
\[
(\kappa m_i)H\left(\frac{|v_i|}{\kappa m_i}\right) \ge
\begin{cases}
\displaystyle
 (\kappa m_{i}) H\left(\frac{\varphi(\kappa)}{\kappa m_{i}}\right) \ge
 (1-\eta) \frac{\varphi(\kappa)^2}{2\kappa m_{i}}, & \mbox{if }|u_i-\kappa m_i|\ge\varphi(\kappa), \\
 0, & \mbox{otherwise},
\end{cases}
\]
where $\eta=\eta(\kappa_0)>0$.

For the second case $u_i=\kappa m_i+v_i=0$, $\log(\kappa m_i+v_i)!/(\kappa m_i)!$ is bounded below by
\[
 v_i \log(\kappa m_i) - \frac{1}{2}\log(\kappa m_i) 
+ \kappa m_i - \frac{1}{2}\log 2\pi - \frac{1}{12}.
\]
Note that the third term is $\kappa m_i=O(\kappa)$, which is of larger order than the corresponding bound for $u_i\ge 1$, i.e.,
$(1-\eta)\varphi(\kappa)^2/(2\kappa m_i)=o(\kappa^{1/3})$.

Because of the assumption that at least one $i$ exists such that $|u_i-\kappa m_i|\ge\varphi(\kappa)$, by summing with respect to $i$, we have
\[
 \log \frac{e^{(u-\kappa m)\cdot\xi} (\kappa m)!}{u!}
\le \sum_{i=1}^n\frac{1}{2}\log(\kappa m_i) - \frac{(1-\eta)\varphi(\kappa)^2}{2\kappa \max m_i} + O(1)
\]
and
\[
 \frac{e^{(u-\kappa m)\cdot\xi} (\kappa m)!}{u!}
=O\left(\kappa^{n/2}\exp\left(-\frac{(1-\eta)\varphi(\kappa)^2}{2\kappa\max m_i}\right)\right).
\]
Since
\[
 \# \{ u\in S_\kappa \,|\, \exists i,\,|u_i-\kappa m_i|\ge\varphi(\kappa) \}
\le \# S_\kappa = O(\kappa^{n-d}),
\]
we have
\begin{equation}
\label{eq:tail_bound}
 \sum_{\exists i,|u_i-\kappa m_i|\ge\varphi(\kappa)}
 \frac{e^{(u-\kappa m)\cdot\xi} (\kappa m)!}{u!}
 =O\left(\kappa^{3n/2-d}\exp\left(-\frac{(1-\eta)\varphi(\kappa)^2}{2\kappa\max m_i}\right)\right).
\end{equation}
From (\ref{eq:riemann}), (\ref{eq:normalizing}), 
and (\ref{eq:tail_bound}) with choosing the function $\varphi(\kappa)$ to be $\kappa/\varphi(\kappa)^2=o(1)$, 
we get an asymptotic evaluation for $P_\kappa(\kappa m,\xi)$.
Therefore, we have
\begin{equation} \label{lm:normalizing}
 P_\kappa(\kappa m,\xi) \sim \frac{\det(\A M^{-1}\A^T)^{1/2}}{(2\pi\kappa)^{(n-d)/2}}
\end{equation}
as $\kappa\to\infty$.
Theorem \ref{th:asymptotic} follows from (\ref{lm:asymptotic1}) and (\ref{lm:normalizing}).
\qed

The approximation of Theorem \ref{th:asymptotic} is interpreted that $U=(U_1,\ldots,U_n)^T$ is distributed as a degenerate normal distribution with the mean vector $m=(m_1,\ldots,m_n)^T$.
Starting from (\ref{lm:asymptotic1}) and applying the approximation arguments of (\ref{lm:normalizing}) again, we can prove that, as $\kappa\to\infty$,
\begin{align*}
 \sum_{u\in S_\kappa} u P_\kappa(u,\xi)
 \sim
 \int_{\{v \,|\, Av=0\}} 
 (\kappa m + v) \widehat P_\kappa(\kappa m+v,\xi) dv 
 \sim
 \kappa 
m(\lambda). 
\end{align*}
\begin{theorem}  \label{th:EisM}
We retain the assumption of Theorem \ref{th:asymptotic}.
For each $\xi\in\R^n$,
\[
 \lim_{\kappa\to\infty} 
 \frac{1}{\kappa}
 \sum_{u\in S_\kappa} u P_\kappa(u,\xi) = 
 m(\lambda),
\]
where $\lambda=\A\xi$ and $m(\lambda)$ is the solution (\ref{eq:ips}) of the IPS.
\end{theorem}
This theorem means that the moment map
$\xi \mapsto \sum_{u\in S_\kappa} u P_\kappa(u,\xi)$
is asymptotically equivalent to the IPS procedure
$\xi \mapsto \kappa m(\A\xi)$.

\begin{example} \rm \label{ex:2F1b}
\comment{The ips is a method to find the MLE.
We note that the output of the standard ips does not agree 
with our expectation in general.}
This is a continuation of Example \ref{ex:2F1a}.
We mean by the IPS the following iteration procedure:
Set
$$ 
 m'_{ij} = m_{ij} \frac{u_{i+}}{m_{i+}}, \ 
 m''_{ij} = m'_{ij} \frac{u_{+j}}{m'_{+j}}
$$
and use $m''$ as the new $m_{ij}$ for the next step.
Here, $u_{i+}$ denotes the $i$-th row sum of $u$ and $u_{+j}$ denotes the $j$-th column sum of $u$.
When the initial $p=m$ satisfies $\frac{m_{12} m_{21}}{m_{11} m_{22}} = 1$,
we can see that the output agrees with 
$E[U](p) = (111/4, 33/4, 37/4, 11/4) \simeq (27.75, 8.25, 9.25, 2.75)$.
However, when the ratio is not equal to $1$, they do not agree in general.
Let $(36,12)$, $(37,11)$ be the row sums and the column sums, respectively.
The expectation of $U$ at $z=\frac{p_{12} p_{21}}{p_{11} p_{22}}=1/2$
is 
\begin{eqnarray*}
& & \left(
   \frac{6595942429}{227713625},\frac{1601748071}{227713625}, 
   \frac{1829461696}{227713625},\frac{903101804}{227713625} \right) \\
& \simeq &
(28.966,7.03405,8.03405,3.96595).
\end{eqnarray*}
Note that it is a vector of rational numbers.
On the other hand, it is known that the $m_{ij}$'s in the steps of the IPS satisfy
the relation
$$ (1/2) m_{11} m_{22}- m_{12} m_{21} = 0
$$
when the initial value of $m$ satisfies $\frac{m_{12}m_{21}}{m_{11} m_{22}} = 1/2$,
and the limit $m_{ij}$'s satisfy
$$
  m_{11}+m_{12} = 36,
 m_{21}+m_{22} = 12,
 m_{11}+m_{21} = 37,
 m_{12}+m_{22} = 11.
$$
(see, e.g., \cite[p.53]{hirotsu} \cite{sinkhorn-knopp}). 
By computing the lexicographic Gr\"obner basis, we find that $m_{11}$ satisfies
the algebraic equation
$m_{11}^2-121 m_{11}+2664=0$,
which does not have a rational solution.
The limit $m$ is approximately equal to
$(28.936572,7.063428,8.063428, 3.936572)$.
Note that it is close to the value of the expectation for $z=1/2$, but differs from it.
When the marginal sums go to infinity, our expectation vector converges
to the IPS value
(see, e.g., \cite{cornfield1956}, \cite[p.21]{hirotsu}, Theorem \ref{th:EisM}).
\end{example}

Noting that (\ref{eq:summation}) can be written  $P(\kappa m,\xi)^{-1} = Z(\kappa\beta;p)\times (\kappa m)!/\exp(\xi\cdot (\kappa m))$, we can obtain an approximate value of the normalizing constant
or the $A$-hypergeometric polynomial by (\ref{lm:normalizing}).
\begin{theorem} \label{th:approxZ}
We retain the assumption of Theorem \ref{th:asymptotic}.
We fix $p$ and $\beta$.
There exists a unique $m \in {\R}_{>0}^n$ such that
$Am =\beta$, $m^{{\bar a}_i} = p^{{\bar a}_i}$ (IPS).
When $\kappa \rightarrow +\infty$, we have
$$ Z(\kappa\beta;p) \sim
 \frac{\left(\prod p_i^{m_i}\right)^\kappa}{\Gamma(\kappa m+1)}
 \frac{(2\pi \kappa)^{n-d}}
      {\det({\bar A} M^{-1} {\bar A}^T)^{1/2}},
$$
where $ M = {\rm diag}(m)$.
\end{theorem}
\begin{example}\rm
 We compare approximate values of $Z$ evaluated by Theorem \ref{th:approxZ}
 with the exact values evaluated by the discrete HGM \cite{goto}, \cite{ott}, \cite{ot2}.
 We consider the $2 \times 4$ contingency tables
 with the fixed row sums $(4, 19)$ and the fixed column sums
 $(9,5,3,6)$. The matrix $A$ is defined
 as in section \ref{sec:classical} (see also Example \ref{ex:2x3}).
 We set $p=(1,1/3,12,1/5001,1,1,1,1)$.
 Note  that $p_4$ is set to be smaller than the other $p_i$'s.
 The constant vector $m$ is determined by the fifth IPS iteration.
 This $m$ is equal to
 $$( 2.79518, 0.652785, 0.551505, 0.000540425,
 6.20482, 4.34722, 2.4485, 5.99946 ).$$
 The exact value of the expectation by the HGM is 
 $$(2.83214,0.627808,0.539555,0.000496547,6.16786,4.37219,2.46044,5.9995).
 $$
 The ratios of the IPS values to the exact values are
 $$(0.98695,1.03978,1.02215,1.08837,1.00599,0.994289,0.995147,0.999993).$$
 The following table illustrates that when $\kappa$ approaches $+\infty$,
 the approximate value converges to the exact value.
\begin{center}
\begin{tabular}{cccc}
$\kappa$ & $\log Z$ & Approx by Th \ref{th:approxZ}  & $|\mbox{error}|$ \\ \hline
\ \ \ 9 & \ \ \ $-568.0127$ & \ \ $-569.8179$  & 1.8052 \\
200 & $-26598.4556$  & $-26598.9446$ & 0.4890 \\
300 & $-42685.5415$  & $-42685.9149$ & 0.3734 \\
\end{tabular}
\end{center}
Note that the approximation of $\log Z$ is close to the exact value,
but the evaluation of the probability by this approximate value has
a relatively big error.
For example, consider the table 
$u=(33,1,1,1,48,44,26,53)$
for the case $\kappa=9$.
The exact probability of getting this $u$ is
$3.26465 \times 10^{-7}$.
The probability evaluated by the approximate value of $\log Z$ above
is 
$1.98529 \times 10^{-6}$.
The approximate value is about 6 times larger than the exact value.
\end{example}

\comment{
\qqq{
Comments on classical moment problems.
Discuss on differences with the following works or with the IPS estimation.} \\
N.Eriksson, S.Fienberg, A.Rinaldo, S.Sullivant,
Polyhedral Conditions for the Nonexistence of the MLE for Hierarchical Log-linear Models, \\
D.Geiger, C.Meek, B.Sturmfels,
On the Toric Algebra of Graphical Models,
The Annals of Statistics, 34 (2006), 1463--1492.
}

\bigbreak
\bigbreak
{\it Acknowledgements}.
The authors are grateful to Satoshi Aoki for a comment on information geometry,
to Hidenao Iwane for a comment on the use of quantifier elimination
to check Theorem \ref{th:fdimage} for small $m$, 
and to Tomohide Terasoma and Keiji Matsumoto for comments on analogous problems
(moment maps in GIT and co-Schwartz maps)
in algebraic geometry.
This work was supported by JSPS KAKENHI Grand Numbers
25287018,
25220001.

\end{document}